\documentstyle[12pt]{article}
\newlength{\abstractwidth}
\renewcommand{\thefootnote}{\fnsymbol{footnote}}
\renewcommand{\thanks}[1]{\footnote{#1}} 
\newcommand{\starttext}{
\setcounter{footnote}{0}
\renewcommand{\thefootnote}{\arabic{footnote}}}

\newcommand{\be}{\begin{equation}}
\newcommand{\bea}{\begin{eqnarray}}
\newcommand{\eea}{\end{eqnarray}}
\newcommand{\ee}{\end{equation}}

\newcommand{\<}{\langle}
\renewcommand{\>}{\rangle}
\def\ba{\begin{eqnarray}}
\def\ea{\end{eqnarray}}

\def\al{\alpha}
\def\b{\beta}

\def\d{\delta}

\def\g{\gamma}

\def\m{\mu}

\def\v{\vskip .1in}

\def\[{{\bf [}}
\def\]{{\bf ]}}

\def\pl{\partial}

\begin{document}
\starttext
\baselineskip=18pt
\setcounter{footnote}{0}

\begin{center}
{\Large \bf ON THE K\"AHLER-RICCI FLOW ON } \\
\bigskip
{\Large \bf COMPLEX SURFACES}
\footnote{Research supported in part by National Science Foundation
grants  DMS-02-45371 and DMS-01-00410}

\bigskip\bigskip

{\large D.H. Phong$^*$ and Jacob Sturm$^{\dagger}$} \\

\bigskip

$^*$ Department of Mathematics\\
Columbia University, New York, NY 10027\\

\v

$^{\dagger}$ Department of Mathematics \\
Rutgers University, Newark, NJ 07102

\end{center}

\baselineskip=15pt
\setcounter{equation}{0}
\setcounter{footnote}{0}

\section{Introduction}
\setcounter{equation}{0}

One of the most important properties of a geometric flow is whether it
preserves
the positivity of various notions of curvature. In the case of the K\"ahler-Ricci flow, the positivity of the curvature operator
(Hamilton \cite{H4}), the positivity of the biholomorphic sectional
curvature (Bando \cite{B}, Mok\cite{M}), and the positivity of the scalar
curvature (Hamilton \cite{H1}) are all preserved. However, whether the
positivity of the Ricci curvature is preserved is still not known. As
stressed for example in
Chen-Tian \cite{CT}, this is central to the problem of convergence of the
K\"ahler-Ricci flow on K\"ahler-Einstein manifolds of positive curvature.
The existence of K\"ahler-Einstein metrics
has been conjectured by S.T.Yau \cite{Y} to be
equivalent to stability in geometric invariant theory, and there is strong interest in relating these notions to the behavior of the K\"ahler-Ricci flow.

\medskip

In this note, we show that the positivity of the Ricci curvature
is preserved on compact complex surfaces, under the additional
assumption that the sum of any two eigenvalues of the traceless
curvature operator on traceless $(1,1)$-forms is non-negative.

\section{The curvature operator in the K\"ahler case}

Let $X$ be an $n$-dimensional compact complex manifold, with
a K\"ahler metric $ds^2=g_{\bar kj}dz^jd\bar z^k$. The K\"ahler-Ricci flow
is
the flow $\dot g_{\bar kj}=-R_{\bar kj}+\mu g_{\bar kj}$, where $n\mu$ is
the
average scalar curvature, and $R_{\bar kj}$ is the Ricci curvature.
By differentiating the defining relation $[\nabla_j,\nabla_{\bar
k}]V^p=R_{\bar
kj}{}^p{}_qV^q$, we obtain the
corresponding flows for the Riemann curvature tensor $R_{\bar
kj}{}^p{}_q$,
the Ricci curvature $R_{\bar kj}=R_{\bar kj}{}^p{}_p$,
and the scalar curvature $R=g^{j\bar k}R_{\bar kj}$:
\v
\bea
&&
\dot R
=\Delta R-\mu R +R^{\bar kj}R_{\bar kj}
\nonumber\\
&&
\dot R_{\bar kj}
=
\Delta R_{\bar kj}
+R_{\bar kl}{}^m{}_jR^l{}_m
-
R_{\bar k}{}^{\bar m}R_{\bar mj}
\nonumber\\
&&
\dot R_{\bar qj\bar lm}
=
\Delta R_{\bar qj\bar lm}+\mu R_{\bar qj\bar lm}
-
R_{\bar l}{}^{\bar r}R_{\bar qj\bar rm}
-
R_{\bar q}{}^{\bar r}R_{\bar lj\bar rm}
+
R_{\bar r}{}^{\bar p}{}_{\bar qj}R_{\bar p}{}^{\bar r}{}_{\bar lm}
\nonumber\\
&&
\quad\quad\quad\quad\quad
+
R_{\bar q}{}^{\bar p}{}_{\bar rm}R_{\bar p}{}^{\bar r}{}_{\bar lj}
-
R_{\bar q}{}^{\bar p}{}_{\bar l}{}^{\bar r}R_{\bar p j\bar r m}
\eea
Here $\Delta=\nabla_l\nabla^l=g^{l\bar k}\nabla_l\nabla_{\bar k}$ is the
complex
Laplacian. It is easily seen that the flows of $R$
and $R_{\bar kj}$ can be written in the same form with $\Delta$
replaced by $\bar\Delta=\nabla_{\bar l}\nabla^{\bar l}=\nabla^l\nabla_l$,
and
hence with ${1\over 2}\Delta_{\bf R}
={1\over 2}(\Delta+\bar\Delta)$. On the other hand, the flow for
the Riemann curvature tensor becomes, when written with $\bar\Delta$
\be
\dot R_{\bar qj\bar lm}
=\bar\Delta R_{\bar qj\bar lm}+\mu R_{\bar qj\bar lm}
-R^r{}_mR_{\bar qj\bar lr}-R^r{}_jR_{\bar qm\bar lr}
+
R^p{}_{r\bar qj}R^r{}_{p\bar lm}+R^p{}_{j\bar lr}R^r{}_{p\bar qm}
-
R^p{}_j{}^r{}_mR_{\bar q p\bar lr}
\ee
Combining the flows with $\Delta$ and $\bar\Delta$, we obtain
the flow with {\it real} Laplacian:
\bea
\dot R_{\bar qj\bar lm}
&=&{1\over 2}\Delta_{\bf R} R_{\bar qj\bar lm}+\mu R_{\bar qj\bar lm}
-{1\over 2}(R^r{}_mR_{\bar qj\bar lr}+R^r{}_jR_{\bar qm\bar lr}
+
R_{\bar l}{}^{\bar r}R_{\bar qj\bar rm}
+
R_{\bar q}{}^{\bar r}R_{\bar lj\bar rm})
\nonumber\\
&&
+
R^p{}_{r\bar qj}R^r{}_{p\bar lm}
+
R^p{}_{r\bar lj}R^r{}_{p\bar qm}
-
R_{\bar q}{}^{\bar p}{}_{\bar l}{}^{\bar r}
R_{\bar pj\bar rm}
\eea

\medskip
As in the Riemannian case \cite{H3,H4}, the flow for the Riemann curvature
operator simplifies considerably in the formalism of frames.
Let $e_a={\pl\over\pl z^j}e^j{}_a$, $e_{\bar a}=e_{\bar a}{}^{\bar
j}{\pl\over\pl \bar z^j}$ be an orthonormal frame at time $t=0$, i.e.,
$e_{\bar b }{}^{\bar k}g_{\bar kj}e^j{}_a=\delta_{\bar b a}$,
$e^j{}_ae^a{}_k=\delta^j{}_k$,
and $e^j{}_a\delta^{a\bar b}e_{\bar b}{}^{\bar k}=g^{j\bar k}$.
Let $g_{\bar kj}$ flow by $\dot g_{\bar kj}$. We want to flow $e^j{}_a$ so
that it remains an orthonormal frame with time. Thus we impose
$0=(g_{\bar kj}e^j{}_ae_{\bar b}{}^{\bar k})\dot{}
=
\dot g_{\bar kj}e^j{}_ae_{\bar b}{}^{\bar k}
+
g_{\bar kj}\dot e^j{}_ae_{\bar b}{}^{\bar k}
+
g_{\bar kj}e^j{}_a\dot e_{\bar b}{}^{\bar k}$.
For this to hold, it suffices to set
\be
\dot e^j{}_a=-{1\over 2}g^{j\bar r}\dot g_{\bar rs}e^s{}_a
=
-{1\over 2}g^{j\bar r}(-R_{\bar r s}+\mu g_{\bar rs})e^s{}_a
\ee
from which it follows that
$\dot e_{\bar b}{}^{\bar k}
={1\over 2}R_{\bar b}{}^{\bar k}-{1\over 2}\mu e^{\bar k}{}_{\bar b}$,
$\dot e^b{}_k=-{1\over 2}R^b{}_k+{1\over 2}\mu e^b{}_k$,
where in general, we can go back and forth between middle Latin indices
($j$,$k$,$l$...) and early
Latin indices ($a$,$b$,$c$...) by using frames, e.g.
$V^a=e^a{}_jV^j$,
$R_{\bar a b\bar cd}=e_{\bar a}{}^{\bar j}e^k{}_{b}e_{\bar c}{}^{\bar l}
e^m{}_{d}R_{\bar jk\bar lm}$.

\medskip

The flow of the frame gets rid of all
the terms mixing the Ricci tensor and the curvature tensor
in the flow of $R_{\bar ab\bar cd}$. Indeed, the cancellation mechanism is
very
simple:
\bea
\dot R_{\bar a b\bar cd}&=&\dot e_{\bar a}{}^{\bar j}e^k{}_{b}e_{\bar
c}{}^{\bar
l}
e^m{}_{d}R_{\bar jk\bar lm}
+
e_{\bar a}{}^{\bar j}\dot e^k{}_{b}e_{\bar c}{}^{\bar l}
e^m{}_{d}R_{\bar jk\bar lm}
+
e_{\bar a}{}^{\bar j}e^k{}_{b}\dot e_{\bar c}{}^{\bar l}
e^m{}_{d}R_{\bar jk\bar lm}
\nonumber\\
&&+
e_{\bar a}{}^{\bar j}e^k{}_{b}e_{\bar c}{}^{\bar l}
\dot e^m{}_{d}R_{\bar jk\bar lm}
+
e_{\bar a}{}^{\bar j}e^k{}_{b}e_{\bar c}{}^{\bar l}
e^m{}_{d}\dot R_{\bar jk\bar lm}
\eea
We have for example
$\dot e_{\bar a}{}^{\bar j}e^k{}_{b}e_{\bar c}{}^{\bar l}
e^m{}_{d}R_{\bar jk\bar lm}
={1\over 2}R_{\bar a}{}^{\bar q}R_{\bar qb\bar cd}-{1\over 2}\mu R_{\bar
ab\bar
cd}$,
and the first term on the right hand side cancels with one of the
terms in the flow with {\it real} Laplacian.
Altogether, we obtain the equation
\be
\dot R_{\bar ab\bar cd}={1\over 2}\Delta_{\bf R}R_{\bar ab\bar cd}-\mu
R_{\bar
ab\bar cd}
+
R_{\bar p r\bar a b}R_{\bar r p\bar c d}
+
R_{\bar a p\bar r d}R_{\bar pb\bar cr}
-
R_{\bar ap\bar cr}R_{\bar pb\bar rd}
\ee
Similarly, the same simplification occurs for the flow of the Ricci curvature, written in a frame. Differentiating the equation
$R_{\bar ab}=e_{\bar a}{}^{\bar q}e^j{}_bR_{\bar qj}$, we see,
not surprisingly, that the term involving the square of the Ricci curvature cancels
\be
\dot R_{\bar ab}={1\over 2}\Delta_{\bf R}R_{\bar ab}
-\mu R_{\bar ab}+R_{\bar ab}{}^p{}_rR^r{}_p.
\ee

\subsection{The traceless curvature operator $S_{\bar ab\bar cd}$}

To analyze the flow of the Riemannian curvature tensor in the operator
case, it
is convenient to separate out the traces. Thus set
\bea\label{definition of S}
S_{\bar ab}&=& R_{\bar ab}-{1\over n}R\delta_{\bar ab}
\nonumber\\
S_{\bar ab\bar cd}&=&
R_{\bar ab\bar cd}-{1\over n}(R_{\bar ab}\delta_{\bar cd}
+R_{\bar cd}\delta_{\bar ab})+{1\over n^2}R\delta_{\bar ab}\delta_{\bar
cd}
\eea
Then $S_{\bar aa}=0$, $S_{\bar aa\bar cd}=0=S_{\bar ab\bar cc}$,
and a straightforward calculation shows that the flows for
$R$, $R_{\bar ab}$, $R_{\bar ab\bar cd}$ are equivalent to the
following flows for $R$, $S_{\bar ab}$, $S_{\bar ab\bar cd}$
\bea
\label{systemS}
&&
\dot R
={1\over 2}\Delta_{\bf R} R +S_{\bar pr}S_{\bar rp}+{1\over n}R(R-\mu n)
\nonumber\\
&&
\dot S_{\bar ab}
=
{1\over 2}\Delta_{\bf R} S_{\bar ab}
+{1\over n}(R-\mu n)S_{\bar ab}
+S_{\bar ab\bar cd}S_{\bar dc}
\\
&&
\dot
S_{\bar ab\bar cd}
=
{1\over 2}\Delta_{\bf R}S_{\bar ab\bar cd}
-\mu S_{\bar ab\bar cd}
+
S_{\bar pr\bar ab}S_{\bar rp\bar cd}
+
S_{\bar ap\bar rd}S_{\bar pb\bar cr}
-
S_{\bar ap\bar cr}S_{\bar pb\bar rd}
+
{1\over n}S_{\bar ab}S_{\bar cd}
\nonumber
\eea
In the K\"ahler case, the Riemann curvature tensor can be viewed as a
symmetric
operator $Op(R)$ on the space $\Lambda^{1,1}$ of real $(1,1)$-forms. This space
itself decomposes into the line spanned by the K\"ahler form
$\omega={\sqrt{-1}\over
2}g_{\bar kj}dz^j\wedge d\bar z^k$, and its orthogonal complement, namely
the
space $\Lambda_0^{1,1}$ of traceless real $(1,1)$-forms. Now
the term $R_{\bar p r\bar a b}R_{\bar r p\bar c d}$ can clearly be viewed as
$Op(R)^2$. Similarly, the tensor $S_{\bar ab\bar cd}$ can be viewed as an
operator $Op(S)$ on $\Lambda_0^{1,1}$,
and we have the decomposition
\be
Op(R)=\pmatrix{R/2 & S\cr S^t& Op(S)\cr}
\ee
The term $S_{\bar pr\bar ab}S_{\bar rp\bar cd}$ in the flow for $S_{\bar ab\bar cd}$ corresponds to $Op(S)^2$.
Following Hamilton \cite{H3,H4}, we show
that the remaining terms $S_{\bar a p\bar r d}S_{\bar pb\bar cr}-S_{\bar
ap\bar
cr}S_{\bar pb\bar rd}$ admit a Lie algebra interpretation. Define the Lie
bracket by
\be
[\phi,\psi]_{\bar ab}
=
\phi_{\bar a p}\psi_{\bar pb}-
\psi_{\bar ap}\phi_{\bar pb}
\ee
Let $\phi_{\bar ab}^{\alpha}$ be an orthonormal basis of real traceless
$(1,1)$-forms, and set $S_{\bar ab\bar cd}=\sum_{\alpha\beta}
M_{\alpha\beta}\phi_{\bar ab}^{\alpha}\phi_{\bar cd}^{\beta}$. Thus
$M_{\alpha\beta}$ is the matrix of $Op(S)$ in the basis
$\phi_{\bar ab}^{\alpha}$. Then
\be
S_{\bar ap\bar cr}S_{\bar pb\bar rd}-S_{\bar a p\bar r d}S_{\bar pb\bar
cr}
=
M_{\al\lambda}M_{\beta\mu}\phi_{\bar ap}^{\alpha}\phi_{\bar pb}^{\beta}
(\phi_{\bar cr}^{\lambda}\phi_{\bar rd}^{\mu}
-
\phi_{\bar rd}^{\lambda}\phi_{\bar cr}^{\mu})
=
M_{\al\lambda}M_{\beta\mu}\phi_{\bar ap}^{\alpha}\phi_{\bar pb}^{\beta}
[\phi^{\lambda},\phi^{\mu}]_{\bar cd}
\ee
Set $[\phi^{\lambda},\phi^{\mu}]=c^{\lambda\mu\rho}\phi^{\rho}$,
where $c^{\lambda\mu\rho}$ are the structure constants of the Lie algebra. The
antisymmetry of $c^{\lambda\rho\mu}$ implies
$M_{\al\lambda}M_{\beta\mu}\phi_{\bar ap}^{\alpha}\phi_{\bar pb}^{\beta}
c^{\lambda\mu\rho}\phi_{\bar cd}^{\rho}=
{1\over 2}M_{\al\lambda}M_{\beta\mu}[\phi^{\alpha},\phi^{\beta}]_{\bar
ab}c^{\lambda\mu\rho}\phi_{\bar cd}^{\rho}$, and thus
\be
S_{\bar ap\bar cr}S_{\bar pb\bar rd}-S_{\bar a p\bar r d}S_{\bar pb\bar
cr}
=
{1\over 2}M_{\al\lambda}M_{\beta\mu}c^{\alpha\beta\nu}\phi_{\bar ab}^{\nu}
c^{\lambda\mu\rho}\phi_{\bar cd}^{\rho}
\equiv {1\over 2}M_{\nu\rho}^{\#}\phi_{\bar ab}^{\nu}\phi_{\bar cd}^{\rho}
\ee
To make $M_{\nu\rho}^{\#}$ explicit, we need the structure
constants $c^{\alpha\beta\nu}$ of the Lie algebra
of traceless $(1,1)$-forms.
Choose a coordinate system centered at a point $p\in M$ such that
the metric $g_{\bar kj}$ is the identity matrix at $p$. Then an orthogonal
basis
for the space of real $(1,1)$ forms is
(in dimension $2$ to simplify notations)
$\omega=dx_1\wedge dy_1+dx_2\wedge dy_2=
{\sqrt{-1}\over 2}(dz_1\wedge d\bar z_1+dz_2\wedge d\bar z_2)$,
$ \eta_1=dx_1\wedge dy_1-dx_2\wedge dy_2\ = \
{\sqrt{-1}\over 2}(dz_1\wedge d\bar z_1-dz_2\wedge d\bar z_2)$,
$\eta_2= dx_1\wedge dy_2+dx_2\wedge dy_1=
{\sqrt{-1}\over 2}(dz_1\wedge d\bar z_2+dz_2\wedge d\bar z_1)$,
$\eta_3=dx_1\wedge dx_2+dy_1\wedge dy_2={1\over 2}(dz_1\wedge d\bar
z_2-dz_2\wedge d\bar z_1)$, with ${\sqrt 2}\eta_i$ forming an orthonormal
basis for $\Lambda_0^{1,1}$. Furthermore,
$[\eta_2,\eta_3]=\eta_1$,
$[\phi_1,\phi_2]=\phi_3$,
$[\phi_3,\phi_1]=\phi_2$,
which means that $\Lambda_0^{1,1}$ is $su(2)$ with structure constants
\be
c^{\al\beta\gamma}=\sqrt 2\epsilon^{\al\beta\gamma}
\ee
where $\epsilon^{\al\beta\gamma}$ is the sign of the permutation
$(1,2,3)\mapsto (\al,\b,\g) $.

\subsection{Positivity of the Ricci curvature in dimension 2}

We are now in position to prove the following theorem:

\bigskip

{\bf Theorem.}
{\it Let $X$ be a compact K\"ahler manifold of dimension $2$, and consider the
K\"ahler-Ricci flow $\dot g_{\bar kj}=-R_{\bar kj}+\mu g_{\bar kj}$. If the initial metric has Ricci
curvature non-negative everywhere and positive somewhere, and if the sum of the two lowest eigenvalues of the
operator $S_{\bar ab\bar cd}$ on the space $\Lambda_0^{1,1}$ of
traceless $(1,1)$-forms is non-negative, then both of these properties continue to hold for all time $t>0$.}

\bigskip
{\it Proof}. If we view the Ricci curvature as a Hermitian form on
$T^{1,0}$
vectors, its positivity is equivalent to the positivity of its trace and
of its
determinant. Set
\be
R_{\bar ab} \ = \ A{\omega\over \sqrt {-1}/2}+
B_1{\eta_1\over \sqrt {-1}/2}+B_2{\eta_2\over \sqrt {-1}/2}+
B_3{\eta_3\over \sqrt {-1}/2}
\ee
In particular, $A={1\over 2}R$ and ${\sqrt 2\over\sqrt{-1}} B_i$ are the components of $S_{\bar
ab}$ in the orthonormal basis ${\sqrt 2}\eta_i$
for $\Lambda_0^{1,1}$.
\v
Claim: The Ricci curvature is non-negative if and only if $A\geq 0$ and
\be
A^2-B_1^2-B_2^2-B_3^2\geq 0\ ,\ \ \ \hbox{i.e.,  $S_{\bar ab}S_{\bar
ba}\leq
{1\over 2}R^2$}
\ee
\v
To see this, we let $X=a{\pl\over \pl z_1}+b{\pl\over \pl z_2}$
be an arbitrary tangent vector. Then
\bea
Ricci(X,\bar X)&=& A(|a|^2+|b|^2)+B_1(|a|^2-|b|^2)+B_2(a\bar
b+b\bar a)-\sqrt{-1}B_3(a\bar b-b\bar a)
\nonumber\\
&=&\left(
\matrix{a& b
}
\right)
 \left(
\matrix{A+B_1 & B_2-\sqrt{-1} B_3\cr
       B_2+\sqrt{-1}B_3 & A-B_1
}
\right)
 \left(
\matrix{\bar a \cr \bar b
}
\right)=
\left(
\matrix{a& b
}
\right)P \left(
\matrix{\bar a \cr \bar b
}
\right)
\nonumber
\eea
Thus the Ricci curvature is non-negative if and only if the matrix $P$ is non-negative.
Now the trace of $P$ is $2A$ and the determinant is $A^2-B_1^2-B_2^2-B_3^2$.
This proves the claim.

\medskip

Set $|S|^2=S_{\bar pr}S_{\bar rp}$. Using the flow for $S_{\bar pr}$, we find
\be
(|S|^2)\dot{}
=
{1\over 2}\Delta_{\bf R}|S|^2
-
(\nabla^lS_{\bar pr}\nabla_l S_{\bar rp}
+
\nabla_lS_{\bar pr}\nabla^l S_{\bar rp})
+
{2\over n}(R-\mu n)|S|^2
+
2\, S_{\bar ba}S_{\bar ab\bar cd}S_{\bar dc}
\ee
Combining with the flow for $R$
\be
\dot R={1\over 2}\Delta_{\bf R}R+|S|^2+{1\over n}R(R-\mu n),
\ee
we obtain the flow for the determinant of the Ricci curvature
\bea
({1\over 2}R^2-|S|^2)\dot{}
&=&
{1\over 2}\Delta_{\bf R}({1\over 2}R^2-|S|^2)
-
\nabla_lR\nabla^lR
+
(\nabla^lS_{\bar pr}\nabla_l S_{\bar rp}
+
\nabla_lS_{\bar pr}\nabla^l S_{\bar rp})
\nonumber\\
&&
\quad
+
{2\over n}(R-\mu n)({1\over 2}R^2-|S|^2)
+
R|S|^2-2\, S_{\bar ba}S_{\bar ab\bar cd}S_{\bar dc}
\eea
We shall abbreviate this equation by
\bea
({1\over 2}R^2-|S|^2)\dot{}
&=&
{1\over 2}\Delta_{\bf R}({1\over 2}R^2-|S|^2)
-
|\nabla R|^2
+
(\nabla S\bar\nabla\bar S
+
\bar\nabla S\nabla\bar S)
\nonumber\\
&&
\quad
+
{2\over n}(R-\mu n)({1\over 2}R^2-|S|^2)
+R|S|^2-2\< S\,Op(S)\,S\>
\eea

\medskip
We examine the non-negativity of the expression ${1\over 2}R^2-|S|^2$, assuming that
it is non-negative at initial time. Consider then the first time
when ${\rm min}\,({1\over 2}R^2-|S|^2)=0$, and consider a minimum point. At this point, by the maximum principle, we have
\be
({1\over 2}R^2-|S|^2)\dot{}
\geq
-
|\nabla R|^2
+
(\nabla S\bar\nabla\bar S
+
\bar\nabla S\nabla\bar S)
+R|S|^2-2\< S\,Op(S)\,S\>
\ee
On the other hand, at a minimum, the derivatives of ${1\over 2}R^2-|S|^2$ all vanish. Thus we have
\be
\nabla_l R={1\over R}(\nabla_l S\cdot\bar S+S\cdot\nabla_l\bar S)
\ee
and hence
\be
|\nabla_lR|
\leq {1\over R}(|\nabla_l S|\cdot |\bar S|+|S|\cdot|\nabla\bar S|)
 =
 {1\over \sqrt 2}
 (|\nabla_l S|+|\nabla _l\bar S|)
 \ee
 since ${1\over 2}R^2-|S|^2=0$. But then
 \be
 \sum_l|\nabla_lR|^2
 \leq {1\over 2}\sum_l(|\nabla_l S|+|\nabla _l\bar S|)^2
 \nonumber\\
 \leq{}\sum_l(|\nabla_l S|^2+|\nabla _l\bar S|^2)
 \ee
 (In the preceding argument, we have assumed that $R>0$,
 which follows from the strong maximum principle if $t>0$.
 If $R=0$ and $t=0$, then we are at a minimum of $R$,
 and $\nabla_lR=0$, so that the above inequality holds trivially). Thus the inequality from the maximum principle reduces to
 \be
 ({1\over 2}R^2-|S|^2)\dot{}
 \geq
 R|S|^2-2\<S\,Op(S)\,S\>
 \ee
 In an orthonormal basis $\phi_{\bar ab}^{\alpha}$
 for the space of traceless $(1,1)$-forms where the operator
 $S_{\bar ab\bar cd}$ is diagonal, with eigenvalues $m_1,m_2,m_3$,
 the preceding inequality can be rewritten as
 \be
 ({1\over 2}R^2-|S|^2)\dot{}
 \geq
 2\sum_{\al=1}^3({1\over 2}R-m_\al)|s_\al|^2
 \ee
 where we have denoted by $s_\al\in {\bf R}$
 the components of $S_{\bar ab}$ in that basis:
 \be
 S_{\bar ab}=\sum_{\al=1}^3s_{\al}\phi_{\bar ab}^{\al}
 \ee
It follows from (\ref{definition of S})
that $S_{\bar a b\bar b a} = {n-1\over n}R={1\over
2}R$ when $n=2$. On the other hand, since
$ S_{\bar a b\bar c d}= \sum_{\al=1}^3
m_\al \phi_{\bar a b}^{\alpha}\phi_{\bar c d}^{\alpha}
$ we obtain
$S_{\bar a b\bar b a} = \sum_{\al=1}^3
m_\al$.
Thus the non-negativity of the determinant of the Ricci curvature will be preserved if we can show that
 \be
 0\leq {1\over 2}R-m_{\alpha}=\sum_{\beta\not=\alpha}m_{\beta},
 \ee
 that is, the sum of any two eigenvalues of $S_{\bar ab\bar cd}$
 is non-negative.

 \medskip
Recall that a symmetric bilinear form is 2-nonnegative if the sum of its two smallest eigenvalues is non-negative.
We have assumed that the traceless curvature operator $S_{\bar ab\bar cd}$ is 2-nonnegative at initial time.
It remains to show that the 2-nonnegativity of the traceless curvature operator is preserved under the K\"ahler-Ricci flow. 
Chen \cite{Chen} has shown that the 2-nonnegativity of the curvature operator $Op(R)$ is preserved by the Ricci flow. Now if the Riemann curvature operator $Op(R)$ is 2-nonnegative, then so is $M=Op(S)$, but the converse does not hold, so we cannot directly quote Chen's result.
\v
First note that if $m_1\leq m_2\leq m_3$ are the eigenvalues of $M$, then
$$ m_1+m_2\ = \ \inf\{M(\phi,\phi)+M(\psi,\psi) : \phi,\psi\in \Lambda_0^{1,1},
|\phi|=|\psi|=1, \phi\perp\psi\}
$$
Moreover, the condition $m_1+m_2\geq 0$ is clearly closed and convex.
The ODE associated to $M$ from the heat flow for the system
(\ref{systemS}) for
$R, S_{\bar ab}, S_{\bar ab\bar cd}$ is
\be 
\label{ODE}{dM\over dt}\ = \ -\m M+M^2+M^\# + T
\ee
where, in coordinates where $M$ is diagonal,
$M_{\al\b}^{\#} = -2(\prod_{\g\not= \al} m_\g)\d_{\al\b}$ and
$T_{\al\b}=s_\al s_\b $. To show that $m_1+m_2\geq 0$ is preserved,
it suffices, by Hamilton's maximum principle for systems, to show that
(\ref{ODE}) preserves this condition. Now Lemma 3.5 of \cite{H4} implies
$$ {d\over dt}(m_1+m_2)\ \geq \ \inf\{{dM\over dt}(\phi,\phi)
+{dM\over dt}(\psi,\psi) \} :
$$
where $\phi,\psi$ range over all $\phi,\psi\in \Lambda_0^{1,1}$ such that $
|\phi|=|\psi|=1, \phi\perp\psi$ and  $ M(\phi,\phi)
+M(\psi,\psi)=m_1+m_2$.
For such $\phi,\psi$, we have $M^2(\phi,\phi)
+M^2(\psi,\psi)=m_1^2+m_2^2$,
$M^{\#}(\phi,\phi)
+M^{\#}(\psi,\psi)= -2m_3(m_1+m_2) $ and $T(\phi,\phi)
+T(\psi,\psi)\geq 0$,
since $T$ is a non-negative operator. Thus
 (\ref{ODE}) implies
\be
 {d\over dt}(m_1+m_2)
 \geq
 -\mu(m_1+m_2)+m_1^2+m_2^2
 -
 2m_3(m_1+m_2).
 \ee
The right hand side is non-negative when $m_1+m_2$ becomes $0$. Thus the non-negativity of
$m_1+m_2$,
and hence of
${1\over 2}R^2-|S|^2$ is preserved under
the flow. Q.E.D.
\v
{\it Remark}: By flowing
$({1\over 2}R^2-|S|^2)^{-1} $, one can show, using a similar argument,
that $({1\over 2}R^2-|S|^2)$ is bounded below by a positive constant, if it is positive everywhere at the initial time and
if the traceless curvature operator is 2-nonnegative.

\v
{\bf Acknowledgements} The authors would like to thank Mu-Tao Wang for his comments and encouragement. They would also like to
thank the Centro di Ricerca Matematica Ennio De Georgi for its warm hospitality when part of this research was carried out.

\end{document}